\documentclass[a4paper,10pt]{article}
\usepackage[english]{babel}
\usepackage[nottoc]{tocbibind}
\usepackage[sorting=nty]{biblatex}
\providecommand{\keywords}[1]
{
  \small	
  \textbf{\textit{Keywords---}} #1
}
\title{Metric dimension and Zagreb indices of essential ideal graph of a finite commutative ring}
\author{Jamsheena P and  Chithra A V \\Department of Mathematics, National Institute of Technology Calicut, \\ Kozhikode, 673601, Kerala, India}

\usepackage{amssymb,amsmath}
\usepackage{amsthm,amsfonts}
\usepackage{cleveref}
\usepackage{tikz}
\usepackage[a4paper,top=46mm,bottom=46mm,left=32mm,right=32mm]{geometry}
\usepackage{arydshln}
\usepackage{mathrsfs}%
\usepackage{multirow}
\newtheorem{theorem}{Theorem}[section]
\newtheorem{definition}[theorem]{Definition}
\newtheorem{lemma}[theorem]{Lemma}

\newtheorem{corollary}[theorem]{Corollary}
\newtheorem{example}[theorem]{Example}
\newtheorem{proposition}[theorem]{Proposition}

\begin{document}

\maketitle

\begin{abstract}
   Let $R$ be a commutative ring with unity. The essential ideal graph $\mathcal{E}_{R}$ of $R$ is a graph whose vertex set consists of all nonzero proper ideals of \textit{R}. Two vertices $\hat{I}$ and $\hat{J}$ are adjacent if and only if $\hat{I}+ \hat{J}$ is an essential ideal. 
   In this paper, we characterize the graph  $\mathcal{E}_{R}$ as having a finite metric dimension.
 Additionally, we identify that the essential ideal graph and annihilating ideal graph of the ring $\mathbb{Z}_{n}$ are isomorphic whenever $n$ is a product of distinct primes. We also estimate the metric dimension of the essential ideal graph of the ring $\mathbb{Z}_{n}$. Furthermore, we determine the topological indices, namely the first and the second Zagreb indices, of  $\mathcal{E}_{\mathbb Z_n}$.
\end{abstract}
\keywords{Essential ideal graph, metric dimension, first and second Zagreb indices}\\
\textit{AMS $2010$ Subject Classification}: $ 05C07, 05C12, 05C25$

\section {Introduction}
Let $\Gamma$ be a simple graph with vertex set  $V(\Gamma)=\{v_1, v_2,\cdots,v_n\}$ and edge set $E(\Gamma)$. If a vertex $u$ is adjacent to a vertex $v$ in $\Gamma$, we write $u\sim v$ in $\Gamma$. The set $N(u)= \{v\in V(\Gamma): v\sim u\ \text{in}\ \Gamma\}$, is called the set of \textit{neighbors} of $u$ and $deg(u)= |N(u)|$ is called the \textit{degree} of a vertex $u$. Also, $N[u]=N(u)\cup \{u\}$.
 The \textit{distance} $d(u,v)$ between two vertices $u$ and $v$ of a connected graph $\Gamma$ is the number of edges in the shortest path between $u$ and $v$.
The  \textit{complete graph} $K_{n}$, is a graph in which any two vertices are adjacent.
A graph $\Gamma$ is a $k-\textit{partite graph}$ if $V(\Gamma)$ can be partitioned into $k$ subsets $V_{1},\ V_{2},\ \cdots,\ V_{k}$ (named partite sets) such that the vertices $u$ and $v$ form an edge in $\Gamma$ if they belong to different partite sets.
If, in addition, there exists an edge between every two vertices belonging to different partite sets, then graph $\Gamma$ can be classified as $\textit{complete k-partite graph}$. The graph denoted as $K_{m,n}$ represents a complete bipartite graph consisting of two sets with sizes $m$ and $n$ respectively.
The \textit{induced subgraph}, $\Gamma[S]$, is formed by taking the subset $S$ of vertices from $\Gamma$, along with all the edges that connect vertices solely within $S$.
The \textit{complement} of a graph $\Gamma$ is denoted by  $\overline{\Gamma}$.
The \textit{join} of two graphs, $\Gamma_{1}$ and $\Gamma_{2}$, represented as $\Gamma_{1}\vee \Gamma_{2}$,  is formed by adding edges between any two vertices $v_1$ and $v_2$, where  $v_1\in \Gamma_{1}$ and  $v_2\in \Gamma_{2}$.

 The concept of metric dimension of a graph was introduced by Slater in \cite{slater1975leaves}, and was called locating sets and locating numbers. 
An equivalent terminology was also introduced by Harary and Melter independently
in \cite{harary1976metric}, and used the term resolving set. Slater described the usefulness of these ideas in long-range aids to navigation. Also, these concepts
have some applications in chemistry for representing chemical compounds \cite{johnson1993structure,johnson1998browsable}, or in problems of pattern
recognition and image processing, some of which involve the use of hierarchical data structures \cite{melter1984metric}. Other applications of this concept to the navigation of robots in networks and other areas appear in \cite{chartrand2000resolvability, hulme1984boolean, khuller1996landmarks}. Hence,
according to its applicability resolving sets has become an interesting and popular topic of investigation in graph theory.\\

Topological Indices play a vital role in mathematical chemistry. They give ideas about structural characteristics with easy identification for a molecule. Hence there are a lot of molecular descriptors called graph invariants. A graph invariant is a number that is invariant under graph isomorphisms in graph theory. The graphical invariant is considered as a structural invariant related to a graph. Since the topological index is constructed as a graphical invariant in molecular graph theory, the computing of topological indices of many graph structures has been an attractive research area for scientists especially chemists and mathematicians for a long time \cite{consonni2009molecular,gutman2012mathematical}.
The first and second Zagreb indices of a graph $\Gamma$ introduced in \cite{gutman1972graph}, and elaborated in \cite{gutman1975graph} are degree-based topological indices defined respectively as follows:\\
 $M_1= \displaystyle\sum_{v\in V(\Gamma)} deg(v)^2$ and 
$ M_2= \displaystyle\sum_{\substack{u\sim v\\ u,v\in V(\Gamma)}}  deg(u)deg(v).$

Let $R$ be a commutative ring with nonzero unity. An element $z\in R$ is said to be a zero divisor of $R$ whenever there exists a nonzero element $w\in R$ such that $zw=0$. An ideal $I$ of a ring $R$ is said to be an \textit{annihilating ideal} of $R$ if there exists a nonzero ideal $J$ of $R$ such as $IJ=0$. An ideal $I$ of a ring $R$ which has a nonzero intersection with every other nonzero ideal of $R$ is called an \textit{essential ideal}.

The study of metric dimension and topological indices of graphs related to various algebraic structures has emerged as a compelling area of research in recent times.  In \cite{pirzada2014locating}, S. Pirzada and R. Raja introduced and investigated the metric dimension of the zero divisor graph of a commutative ring $R$. The results on topological indices of this graph can be seen in \cite{selvakumar2022wiener}. In \cite{banerjee2022spectra, banerjee2023adjacency}, S. Banerjee determined the metric dimension and topological indices like the Wiener index, the first and the second Zagreb index of comaximal graph of the ring $\mathbb{Z}_{n}$. In \cite{aijaz2020annihilating}, M. Aijaz and S. Pirzada computed the metric dimension of annihilating ideal graphs of commutative rings.  The annihilating ideal graph $\mathbb{AIG}(R)$, of a commutative ring $R$, introduced and studied by M. Behboodi and Z. Rakeei in \cite{behboodi2011annihilating}, is a graph in which the vertex set consists of the set of all nonzero annihilating ideals of $R$ and two distinct vertices $\hat{I}$ and $\hat{J}$ are joined by an edge if and only if $\hat{I}\hat{J} =0$. 

Being motivated by these works, in this paper, we study the metric dimension and topological indices of the essential ideal graph of the ring  $\mathbb{Z}_{n}$. The essential ideal graph $\mathcal{E}_{R}$ of a commutative ring $R$, introduced and studied by J. Amjadi in \cite{amjadi2018essential}, is a graph in which the vertex set is the set of all nonzero proper ideals of $R$ and two vertices $\hat{I}$ and $\hat{J}$ are joined by an edge whenever $\hat{I}+\hat{J}$ is an essential ideal. To date, there is no information about the metric dimension and topological indices of the essential ideal graph of $\mathbb{Z}_{n}$ in literature.

This paper has been organized as follows: In Section $2$, we list the results and definitions that are needed for the present study. In Section $3$, we determine the metric dimension of the essential ideal graph of $\mathbb{Z}_{n}$. Also, we prove that the essential ideal graph and annihilating ideal graph of the ring $\mathbb{Z}_{n}$ are equal (up to isomorphism) whenever $n$ is a product of $k\ge 2$ distinct primes. Moreover, we provide an alternate proof to show that the metric dimension of $\mathcal{E}_{\mathbb Z_n}$ is $\le k$ when $n=\prod_{i=1}^k p_i$. In section $4$, we calculate the first and the second Zagreb index of the graph $\mathcal{E}_{\mathbb Z_n}$  for any $n\ge 4$.

Throughout this paper, $\mathbb{Z}_{n}= \mathbb{Z}/n\mathbb{Z}$, where $n\ge 4$ and $n$ is not a prime.


\section{Preliminaries}
In this section, we list some definitions and results that are needed for the present study.
\begin{definition}
    A subset $W$ of vertices of a connected graph $\Gamma$ is said to resolve $\Gamma$, if each vertex of $\Gamma$ is uniquely determined by its vector of distances to the vertices of $W$. In general, for an ordered subset $W= \{w_1, w_2,\cdots, w_k\}$ of vertices of a connected graph $\Gamma$ and a vertex $v \in V(\Gamma)\backslash W$ of $\Gamma$, the metric representation of $v$ with respect to $W$ is the $k-$vector $r(v|W)= (d(v,w_1),d(v,w_2), \cdots d(v, w_k))$. The set $W$ is a resolving set for $\Gamma$ if  $r(v|W)\ne r(u|W)$, for any pair of distinct vertices $u,v \in V(\Gamma)\backslash W $.\\

  The resolving set, the metric representation of a vertex, and the metric dimension of a graph are also called the locating set, locating code of a vertex, and locating number of a graph respectively.
\end{definition}
\begin{definition}
   Let $\Gamma$ be a connected graph with order $n\ge 2$. The metric dimension $dim (\Gamma)$ of $\Gamma$, is defined as $dim(\Gamma)= min\{|W|: W \textit{is a resolving set of}\ \Gamma \}$ and such a set $W$ is the metric basis for $\Gamma$.\\
    For every connected graph $\Gamma$ of order $n\ge 2$, $1\le dim(\Gamma)\le n-1$.
\end{definition}
\begin{definition}
  Let $\Gamma$ be a connected graph with order $n\ge 2$. Two distinct vertices $u$ and $v$ are said to be distance similar if $d(u,x)= d(v,x)$, for all $x \in V(\Gamma)\backslash \{u,v\}$. It can be verified that the distance relation is an equivalence relation on $V(\Gamma)$ and two vertices are distance similar if either $uv \notin E(\Gamma)$ and $N(u)= N(v)$ or $uv \in E(\Gamma)$ and $N[u]= N[v]$.
\end{definition}
\begin{theorem}\label{dmsn Pn&Kn} \cite{chartrand2000resolvability}
  Let $\Gamma$ be a connected graph with order $n\ge 2$ and $W$ be a metric basis for $\Gamma$. Then $dim(\Gamma)= n-1$ if and only if $\Gamma \cong K_n$.
\end{theorem}
\begin{theorem}\cite{pirzada2014locating}\label{Dsimilr conds}
Let $\Gamma$ be a connected graph and $V(\Gamma)$ is partitioned into $k$ distinct distance similar classes $X_1, X_2, \cdots X_k$. Then 
\begin{enumerate}
    \item Any resolving set $W$ contains all but at most one vertex from each $X_i$.
    \item If $t$ is the number of distance similar classes that consist of a single vertex, then  $|V(\Gamma)|-k\le dim(\Gamma) \le V(\Gamma)|-k+t$.
\end{enumerate}
\end{theorem}
\begin{theorem}\label{E_R finite}\cite{amjadi2018essential}
   Let $R$ be a commutative ring with unity. Then, $\mathcal{E}_{R}$ is a finite graph if and only if every vertex of $\mathcal{E}_{R}$ has finite degree.
\end{theorem}
In \cite{p2023structure}, the authors determined the structure of essential ideal graph of the ring  $\mathbb{Z}_{n}$ by defining an equivalence relation on the set $\mathscr{U}$ of nonessential ideals of $\mathbb{Z}_{n}$ as follows:  
 \begin{definition}
     Let $\Xi= \{ 1,2, \cdots, k\}$ be an index set.
    For an ideal $\hat{I}$ of $\mathscr{U}$, define a subset  $\Xi_{\hat{I}}$ of $\Xi$  by,
    $\Xi_{\hat{I}}= \{i: r_i= m_i \ \text{in} \ \hat{I}\}$.
\end{definition}
\begin{definition} \label{eqvlnc reln}
    Let $\hat{I}$ and $\hat{J}$ be any two ideals of $\mathscr{U}$. We define a relation $\preccurlyeq$ on $\mathscr{U}$ by $\hat{I}\preccurlyeq \hat{J}$ if and only if $\Xi_{\hat{I}}= \Xi_{\hat{J}}$.
\end{definition}
Thus,  $\mathscr{U}$ is partitioned into $2^k-2$ equivalent classes, and each equivalent class is denoted by $[\hat{I}]$.
For example, if $n= p_1^2p_2^3p_3p_4p_5^4$ and $\hat{I}=  \langle p_2^3p_4 \rangle$ is the representative ideal then, the corresponding equivalent class $[\hat{I}]$ is the set  $X_{\hat{I}}= \{\langle p_1^{r_1}p_2^3p_4p_5^{r_5} \rangle : 0\le r_1\le 1,\ \text{and}\ 0\le r_5\le 3 \}$.

 \begin{lemma}\label{Chrczn Adcncy}
Let $\hat{K}$ and $\hat{L}$ be two vertices of any two of the $2^k -2$ equivalent classes, say $[\hat{I}]$ and $[\hat{M}]$ respectively. Then $\hat{K}$ and $\hat{L}$ are adjacent in $\mathcal{E}_{\mathbb Z_n}$  if and only if $\Xi_{\hat{I}} \cap \Xi_{\hat{M}}= \phi$.
\end{lemma}
The next theorem can be found in \cite{p2023structure}, which determines the structure of the induced subgraph $\mathcal{E}_{\mathbb Z_n}(\mathscr{U})$. 
The next theorem gives the structure of the induced subgraph $\mathcal{E}_{\mathbb Z_n}(\mathscr{U})$. 
\begin{theorem}\label{G-join1} \cite{p2023structure} Let  $n= p_{1}^{m_1}p_{2}^{m_2} \cdots p_{k}^{m_k}$, where $p_{1}<p_{2}< \cdots < p_{k}$ are primes, $k\ge 2$, and $m_i>1$ for at least one $i$.Then, the induced subgraph $\mathcal{E}_{\mathbb Z_n}(\mathscr{U})$ is the generalized join of certain null graphs given by,
\begin{equation*}
\begin{split}
    \mathcal{E}_{\mathbb Z_n}(\mathscr{U}) = \mathscr{G}&[\mathcal{E}_{\mathbb Z_n}([\langle p_1^{m_1}\rangle]),\cdots, \mathcal{E}_{\mathbb Z_n}([\langle p_k^{m_k}\rangle]), \mathcal{E}_{\mathbb Z_n}([\langle p_1^{m_1}p_2^{m_2}\rangle]),\cdots,  \\
    & \mathcal{E}_{\mathbb Z_n}([\langle p_{k-1}^{m_k-1}p_k^{m_k}\rangle]),\cdots, \mathcal{E}_{\mathbb Z_n}([\langle p_2^{m_2}\cdots p_{k-1}^{m_{k-1}}p_k^{m_k} \rangle])],
\end{split}
\end{equation*}
where $\mathcal{E}_{\mathbb Z_n}([\hat{I}])= \overline{K} _ {\displaystyle \prod_{i\notin \Xi_{\hat{I}}}m_i}$ for the representative ideal $\hat{I} (\text{vertex of}\ \mathscr{G})$ of the equivalent class $[\hat{I}]$.
\end{theorem}
The following theorem determines the structure of $\mathcal{E}_{\mathbb Z_n}$ as the join of a complete graph induced by the essential ideals of $\mathbb Z_n$ and the induced subgraph $   \mathcal{E}_{\mathbb Z_n}(\mathscr{U})$.
\begin{theorem} \label{Complt Strctr of EZn}  \cite{p2023structure} 
   Let  $n= p_{1}^{m_1}p_{2}^{m_2} \cdots p_{k}^{m_k}$, where $p_{1}<p_{2}< \cdots < p_{k}$ are primes, and $m_i>1$ for at least one $i$. Then, the essential ideal graph  $\mathcal{E}_{\mathbb Z_n}\cong K_m \vee H$, where $K_m$ is the complete graph  on $m=\prod_{i=1}^{k}m_{i}-1$ vertices and  \begin{equation*}
       \begin{split}
      H= \mathscr{G}&[\mathcal{E}_{\mathbb Z_n}([\langle p_1^{m_1}\rangle]),\cdots, \mathcal{E}_{\mathbb Z_n}([\langle p_k^{m_k}\rangle]), \mathcal{E}_{\mathbb Z_n}([\langle p_1^{m_1}p_2^{m_2}\rangle]),\cdots,  \\
    & \mathcal{E}_{\mathbb Z_n}([\langle p_{k-1}^{m_k-1}p_k^{m_k}\rangle]),\cdots,\mathcal{E}_{\mathbb Z_n}([\langle p_2^{m_2}\cdots p_{k-1}^{m_{k-1}}p_k^{m_k} \rangle])].
\end{split}
\end{equation*}
\end{theorem}

\section{Metric Dimension of $\mathcal{E}_{\mathbb{Z}_{n}}$}
In this section, we compute the metric dimension of the essential ideal graph of $\mathbb{Z}_{n}$.
\begin{theorem}\label{dmn finite R}
Let $R$ be a commutative ring with unity. Then, $dim(\mathcal{E}_{R})$ is finite if and only if $R$ is finite.
\end{theorem}
\begin{proof} If $R$ is finite, obviously, $dim(\mathcal{E}_{R})$ is finite. Conversely, suppose that $dim(\mathcal{E}_{R})=k<\infty$. This ensures that each vertex of $\mathcal{E}_{R})$ has a unique $k$-vector metric representation with respect to a minimum resolving set $W$ of cardinality $k$. Since $diam(\mathcal{E}_{R})=3<\infty$, for every vertex $v\in V(\mathcal{E}_{R})\backslash W$, there are only $4^k$ choices for $r(v|W)$. Hence, $|V(\mathcal{E}_{R})|\le 4^k+k$. \\
\end{proof}
The next result follows directly from Theorems \ref{E_R finite} and \ref{dmn finite R}.
\begin{corollary}
   Let $R$ be a commutative ring with unity. Then, $dim(\mathcal{E}_{R})$ is finite if and only if every vertex of $\mathcal{E}_{R}$ has finite degree.  
\end{corollary}
    \begin{lemma}\label{gcd 1}
        Let $n= p_1p_2\cdots p_k$, where $p_1<p_2<\cdots < p_k$ are primes and  let $d_1$ and $d_2$ be two distinct nontrivial proper divisors of $n$. Then, $gcd(d_1, d_2)=1$ if and only if $n\mid (\frac{n}{d_1})(\frac{n}{d_2})$ 
    \end{lemma}
    \begin{proof}
        Assume that $gcd(d_1,d_2)=1$. Then, there exist integers $x$ and $y$ such that $1= d_1x+d_2y$. Now, \begin{equation*}
        \begin{split}
             (\frac{n}{d_1})=& nx+(\frac{n}{d_1})d_2y\\
            (\frac{n}{d_2})=& (\frac{n}{d_2})d_1x+d_2y \\
            \text{and hence}\ (\frac{n}{d_1})(\frac{n}{d_2})=& n(x_1+ 2nxy+y_1),\\         \end{split}
    \end{equation*} where $x_1= (\frac{n}{d_2})d_1 x^2, \ y_1= (\frac{n}{d_1})d_2y^2$. Thus, $n| (\frac{n}{d_1})(\frac{n}{d_2})$.
        For the converse, suppose that $gcd(d_1,d_2)=d>1$. Then $d=p_{i_1}p_{i_2}\cdots p_{i_t}$, where $p_{i_1}, p_{i_2},\cdots,p_{i_t}$ are primes such that $i_1< i_2< \cdots <i_t$ and $1\le i_t\le k-1$ so that $d_1= r_1 d$, and $d_2=r_2d$. Consequently, both divisors $(\frac{n}{d_1})$ and $(\frac{n}{d_2})$ of $n$ do not have $d$ as a factor and hence  $n\nmid (\frac{n}{d_1})(\frac{n}{d_2})$.
    \end{proof}
    \begin{theorem}\label{isomphsm}
        Let $R_1= \displaystyle \prod_{i=1}^{k}F_i$, where each $F_i$ is a field and let $R_2=\mathbb{Z}_{n}$ for $n=p_1p_2\cdots p_k$, where $p_i$'s are distinct primes for $1\le i \le k$. Then, $\mathcal{E}_{R_{1}}\cong \mathcal{E}_{R_{2}}\cong \mathbb{AIG}(R_2)$.
    \end{theorem}
    \begin{proof}
        We first note that the vertices of $\mathcal{E}_{R_1}$ are the nonzero proper ideals of the ring $\prod_{i=1}^{k}F_i$, given by $\hat{I} =\prod_{i=1}^{k}\hat{I_i}$, where $\hat{I_i}=\langle 0 \rangle$ for at least one $i$ and $\hat{I_i}=F_i$ for at least one $i$. Thus, $|V(\mathcal{E}_{R_1})|= 2^k-2= |V(\mathcal{E}_{R_{2}})|= |V(\mathbb{AIG}(R_2))|$. Also,\\ 
        $ V(\mathcal{E}_{R_{2}})=V(\mathbb{AIG}(R_2))= \{\langle d \rangle: d\ \text{is a positive proper divisor of}\ n\}$. For the divisor $d$ of $n$, define a map $\varphi:V(\mathcal{E}_{R_{2}})\rightarrow V(\mathbb{AIG}(R_2))$ by $d \longmapsto \frac{n}{d}$(divisor conjugate of $d$).
        Since each divisor $d$ of $n$ has a unique divisor conjugate $\frac{n}{d}$, and $1< \frac{n}{d}< n$ for  $1< d< n$,  it follows immediately that   $\varphi$ is both one-one and onto. Now, Lemma \ref{gcd 1} assures that two vertices $\langle d_1\rangle$ and $\langle d_2 \rangle$  are adjacent in $\mathcal{E}_{R_{2}}$ if and only if $\varphi(\langle d_1\rangle)$ and  $\varphi(\langle d_2\rangle)$ are adjacent in $\mathbb{AIG}(R_2)$.  Thus,  $\varphi$ is an isomorphism and hence  $\mathcal{E}_{R_{2}}\cong \mathbb{AIG}(R_2)$.

        Now, for each vertex $\hat{I} =\prod_{i=1}^{k}\hat{I_i}$ of $\mathcal{E}_{R_1}$, we define a subset $\Theta_{\hat{I}}$ of the index set $\{1,2,\cdots,k\}$ such that  $\hat{I_i}=\begin{cases}
            \langle 0 \rangle, & \text{if}\ i\in \Theta_{\hat{I}} \\
            F_i, & \text{otherwise}. 
        \end{cases}$ Obviously, two distinct vertices $\hat{I}$ and $\hat{J}$ are adjacent in $\mathcal{E}_{R_1}$ if and only if $\Theta_{\hat{I}}\cap \Theta_{\hat{J}}=\phi$. Define a map $\psi: V(\mathcal{E}_{R_1})\rightarrow V(\mathbb{AIG}(R_2))$ by $\psi(\hat{I})= \langle \displaystyle\prod_{i\notin \Theta_{\hat{I}}} p_i \rangle$. Clearly, $\psi$ is a well defined bijection preserving adjacencies and nonadjacencies in  $\mathcal{E}_{R_1}$ and $\mathbb{AIG}(R_2)$, and hence $\mathcal{E}_{R_{1}}\cong \mathbb{AIG}(R_2)$.
    \end{proof}
In \cite{aijaz2020annihilating}, the authors computed the metric dimension of the annihilating ideal graph of the rings $\displaystyle \prod_{i=1}^{k}F_i$ and $\mathbb{Z}_{n}$, $n=p_1p_2\cdots p_k$.
    \begin{theorem}\cite{aijaz2020annihilating}\label{dim AIGof Zn}
        For $R=\displaystyle \prod_{i=1}^{k}F_i\ \text{or}\ \mathbb{Z}_{n}$, $n=p_1p_2\cdots p_k$, the following hold:
        \begin{enumerate}
            \item $dim(\mathbb{AIG}(R))=k-1$ for $1\le k \le 4$.
            \item $dim(\mathbb{AIG}(R))=5$ for $k=5$.
            \item $dim (\mathbb{AIG}(R))\le k$ for $k\ge 6$.
        \end{enumerate}
    \end{theorem}
    The proof is developed by showing that the annihilating ideal graph $\mathbb{AIG}(R)$ for $R= \displaystyle \prod_{i=1}^{k}F_i\ \text{or}\ \mathbb{Z}_{n}$, $n=p_1p_2\cdots p_k$ is isomorphic to the zero divisor graph ($\mathbb{ZDG}$) of the boolean ring $\displaystyle\prod_{i=1}^k \mathbb{Z}_{2}$ and applying the result on metric dimension of zero divisor graph of $\displaystyle\prod_{i=1}^k \mathbb{Z}_{2}$ [Proposition $6.2$ and Theorem $6.3$ of \cite{raja2016locating}]. Hence by Theorems \ref{isomphsm} and \ref{dim AIGof Zn}, we can have the following result.

\begin{proposition}
    Let $R=\displaystyle \prod_{i=1}^{k}F_i\ \text{or}\ \mathbb{Z}_{n}$, $n= p_1p_2\cdots p_k$. Then, 

\begin{enumerate}
    \item  $dim(\mathcal{E}_{R})=k-1$ for $1\le k \le 4$.
    \item  $dim(\mathcal{E}_{R})=5$ for $k=5$.
     \item  $dim(\mathcal{E}_{R})\le k$ for $k\ge 6$.
\end{enumerate} 
\end{proposition}
In the following theorem, we give another proof for computing the metric dimension of $\mathcal{E}_{\mathbb{Z}_{n}}$ for $n= p_1p_2\cdots p_k$, $k\ge 6$, by finding a minimal resolving set of $\mathcal{E}_{\mathbb{Z}_{n}}$. For this, we make use of the following Lemma.

\begin{lemma}\label{dist.in EZn}
   Let $R= \mathbb{Z}_{n}$, $n=  p_1p_2\cdots p_k$. Then, for any two vertices $\hat{I}$ and $\hat{J}$ of $\mathcal{E}_{R}$ 
 \begin{enumerate}
  \item $d(\hat{I}, \hat{J})=2$ if and only if $\hat{I}+\hat{J} \ne R$ and $\hat{I}\cap \hat{J}\ne 0$.
 \item $d(\hat{I}, \hat{J})=3$ if and only if  $\hat{I}+\hat{J} \ne R$ and $\hat{I}\cap \hat{J}=0$.
\end{enumerate}
\end{lemma}
\begin{proof}
   $(i)$  First, suppose that $d(\hat{I}, \hat{J})=2$. Obviously, $\hat{I}+\hat{J} \ne R$.  Thus, it remains to prove that $\hat{I}\cap \hat{J}\ne 0$. If possible, let $\hat{I}\cap \hat{J}=0$. Then, any prime not in the generator of the ideal $\hat{I}$ must be in the generator of the ideal $\hat{J}$ and vice versa.
    Hence, if $\hat{K}$ is a vertex adjacent to the vertex $\hat{I}$, then it cannot be adjacent to the vertex $\hat{J}$ as the generators of both $\hat{K}$ and  $\hat{J}$ have at least one common prime factor. This leads to the conclusion that $d(\hat{I}, \hat{J})>2$, is a contradiction. Thus, $\hat{I}\cap \hat{J}\ne 0$. For the  converse, assume that $\hat{I}+\hat{J} \ne R$ and $\hat{I}\cap \hat{J}\ne 0$. Then  $d(\hat{I}, \hat{J})>1$. 
    Since  $\hat{I}\cap \hat{J}\ne 0$, there must exist at least one prime number $p_s$ such that $p_s$ is not a prime factor of generators of both ideals $\hat{I}$ and $\hat{J}$. Thus, if $\hat{S}=\langle p_s\rangle$, we have $\hat{I}\sim \hat{S} \sim \hat{J}$. Consequently,  $d(\hat{I}, \hat{J})= 2$.\\
   $(ii)$ Result follows as a direct consequence of the proof of Case $1$.
\end{proof}
In the following theorem, we give another proof for computing the metric dimension of $\mathcal{E}_{\mathbb{Z}_{n}}$ for $n= p_1p_2\cdots p_k$.
\begin{theorem}
Let $n= p_1p_2\cdots p_k$, where $p_{1}< p_{2}< \cdots < p_{k}$ are primes and $k\ge 6$. Then $dim (\mathcal{E}_{\mathbb{Z}_{n}})\le k$.
\end{theorem}
\begin{proof}
Consider the set $W$ consisting of all  minimal ideals of $\mathcal{E}_{\mathbb{Z}_{n}}$ as in the following order:
   \[ W= \{\langle p_1p_2\cdots p_{k-1}\rangle, \langle  p_1p_2\cdots p_{k-2}p_k\rangle, \cdots, \langle p_2p_3\cdots p_{k}\rangle \}. \] 
  \textbf{claim}: $W$ is a resolving set of $\mathcal{E}_{\mathbb{Z}_{n}}$\\
We need to show that each vertex $v \in V(\mathcal{E}_{\mathbb{Z}_{n}})\backslash W$ has a unique representation of distances with respect to $W$. For this, take any two vertices of $V(\mathcal{E}_{\mathbb{Z}_{n}})\backslash W$ of the form $\hat{I}=\langle p_{i_1}p_{i_2}\cdots p_{i_t}\rangle$ and $\hat{J}=\langle p_{j_1}p_{j_2}\cdots p_{j_s}\rangle$, where $p_{i_1}, p_{i_2},\cdots,p_{i_t},p_{j_1},p_{j_2},\cdots,p_{j_s}$ are primes such that $i_1< i_2< \cdots <i_t$ and $j_1<j_2<\cdots<j_s$ not necessarily distinct and $1\le i_t, j_s \le k-2$. Then three cases may occur- either $t<s$, or $t=s$, or $t>s$. \\
 Case $1:  t<s$\\
Then, there exists at least one prime $p_{j_l}$ which is in the generator of $\hat{J}$ but not in that of $\hat{I}$. Now, consider a vertex $\hat{P}$ in $W$ such that $p_{j_l}$ is not in the generator of $\hat{P}$. Then $d(\hat{I}, \hat{P})=2$, by Lemma \ref{dist.in EZn}$(1)$.
That is, $\hat{I} \sim \langle p_{j_l}\rangle \sim \hat{P}$. However, since $\hat{J}$ is not adjacent to the vertex $\langle p_{j_l}\rangle$ to which $\hat{P}$ is only adjacent, $d(\hat{J}, \hat{P})=3$.  
Then, the coordinate corresponding to the vertex $\hat{P}$ of $W$ in the $k$-vector of both $\hat{I}$ and $\hat{J}$ are distinct. Hence, $r(\hat{I}|W) \ne r(\hat{J}| W)$. \\
Case $2: t=s$ \\
In this case, at least one prime is not common in the generators of both $\hat{I}$ and $\hat{J}$. Without loss of generality, assume that $p_{i_h}$ is in the generator of  $\hat{I}$ but not in that of $\hat{J}$ and $p_{j_l}$ is in the generator of  $\hat{J}$  but not in that of $\hat{I}$. Consider the vertex $\hat{Q}\in W$ such that the generator of  $\hat{Q}$ contains  $p_{j_l}$ as a factor but not $p_{i_h}$. Then, $\hat{Q}$ is adjacent only to the vertex $\langle p_{i_h}\rangle$ and the latter is not adjacent to $\hat{I}$. Hence, by Lemma \ref{dist.in EZn},  $d(\hat{I}, \hat{Q})=3$ and  $d(\hat{J},\hat{Q})=2$.\\
Case $3:  t>s$ \\
Here, there is at least one prime $p_{i_h}$ in the generator of $\hat{I}$ but not in that of $\hat{J}$. Then, by Lemma \ref{dist.in EZn}, $d(\hat{I}, \hat{K})=3$, and $d(\hat{J},\hat{K})= 2$, for the vertex $\hat{K} \in W$ having $p_{i_h}$ not in the generator of $\hat{K}$.
 This proves that $r(\hat{I}|W)\ne r(\hat{J}|W)$, for any two distinct vertices $\hat{I}$ and $\hat{J}$ in $V(\mathcal{E}_{\mathbb{Z}_{n}})\backslash W$. Hence $W$ is a resolving set of cardinality $k$ and $dim (\mathcal{E}_{\mathbb{Z}_{n}})\le k$.
\end{proof} 
\begin{proposition}Let $T=|V(\mathcal{E}_{\mathbb{Z}_{n}})|$. Then, 
    $dim (\mathcal{E}_{\mathbb{Z}_{n}})=T-1$ if and only if either $n= p^m, \ m>1$ or $n=p_1p_2$.
\end{proposition}
\begin{proof}
It is obvious that $dim (\mathcal{E}_{\mathbb{Z}_{n}})=T-1$ when $n= p^m, \ m>1$ or $n=p_1p_2$. For the converse, assume that  $dim (\mathcal{E}_{\mathbb{Z}_{n}})=T-1$. Then, 
$\mathcal{E}_{\mathbb{Z}_{n}}$ is complete by Theorem \ref{dmsn Pn&Kn}. Suppose $n\ne p_1p_2$. To prove  $n= p^m, \ m>1$, assume to the contrary that $n= p_1^{\alpha_1}p_2^{\alpha_2}\cdots p_k^{\alpha_k}$, $k\ge 2$ and $\alpha_i>1$ for at least two $i$ (say, $\alpha_1, \alpha_2$).
 Now, consider the two vertices $\hat{I}=\langle  p_1^{\alpha_1} \rangle$ and $\hat{J}= \langle p_1^{\alpha_1}p_2^{\alpha_2}\rangle$.  Obviously, $\hat{I}$ and $\hat{J}$ are nonadjacent in $\mathcal{E}_{\mathbb{Z}_{n}}$ contradicting the fact that $\mathcal{E}_{\mathbb{Z}_{n}}$ is complete.
\end{proof}

By Theorem \ref{Complt Strctr of EZn}, $\mathcal{E}_{\mathbb{Z}_{n}}\cong K_m \vee \mathscr{G}[\Gamma_1,\Gamma_2,\cdots,\Gamma_{2^k-2}]$, 
 where $\Gamma_i= \mathcal{E}_{\mathbb{Z}_{n}}([\hat{I}])$ for each of the equivalence class $[\hat{I}]$ of the partition on the set of nonessential ideals of $\mathcal{E}_{\mathbb{Z}_{n}}$. This can be further viewed as $\mathcal{E}_{\mathbb{Z}_{n}}\cong  \mathscr{G}[K_m, \Gamma_1,\Gamma_2,\cdots,\Gamma_{2^k-2}]$, since the vertices of the subgraph $K_m$ are adjacent to all the vertices of the subgraphs $\Gamma_i$ for $1\le i \le 2^k-2$. Also, note that the vertices in each of the induced subgraphs $K_m$ and $\Gamma_i$ for $1\le i \le 2^k-2$ are distance similar so that $V( \mathcal{E}_{\mathbb{Z}_{n}})$ is partitioned into $2^k-1$ distance similar classes $X, X_1, X_2, \cdots, X_{2^k-2}$ as follows.
 
 \begin{align*}
    X =& \{\langle p_{1}^{r_1}p_{2}^{r_2} \cdots p_{k}^{r_k}\rangle: 0\le r_i\le m_i-1 \ \text{for}\ 1\le i\le k\}\backslash \mathbb{Z}_{n}, \\
 X_1= & X_{\langle p_{1}^{m_1}\rangle} =  \{\langle p_{1}^{m_1}p_{2}^{r_2} \cdots p_{k}^{r_k}\rangle: 0\le r_i\le m_i-1 \ \text{for}\ 2\le i\le k \},\\
    &\vdots \\
 X_k =& X_{\langle p_{1}^{m_k}\rangle} =  \{\langle p_{1}^{r_1}p_{2}^{r_2} \cdots p_{k}^{m_k}\rangle: 0\le r_i\le m_i-1 \ \text{for}\ 1\le i\le k-1 \},\\
 X_{k+1} =&  X_{\langle p_{1}^{m_1}p_{2}^{m_2}\rangle} =   \{\langle p_{1}^{m_1}p_{2}^{m_2} \cdots p_{k}^{r_k}\rangle: 0\le r_i\le m_i-1 \ \text{for}\ 3\le i\le k \},\\ 
 & \vdots \\
 X_{2^k-2}= & X_{\langle p_{2}^{m_2} \cdots p_{k}^{m_k}\rangle}= \{\langle p_{1}^{r_1}p_{2}^{m_2} \cdots p_{k}^{m_k}\rangle: 0\le r_1\le m_1-1 \}.  
 \end{align*}

 Here, $X= V(K_m)$ and $X_i=V(\Gamma_i)= [\hat{I}]$ for each of $2^k-2$ equivalent class $[\hat{I}]$. By Theorem \ref{Dsimilr conds}, any resolving set $W$ of $\mathcal{E}_{\mathbb{Z}_{n}}$ must contain all but at most one vertex from each of the partitioned sets $X,\ X_i$ for $i= 1,2,\cdots,2^k-2$. Hence, for any resolving set $W$, 
 \begin{equation}
 \begin{split}
      |W|\ge & |X|-1+|X_1|-1+|X_2|-1+\cdots +|X_{2^k-2}|-1 \\
      \ge & m-1+ T-m\underbrace {-1-1\cdots -1}_{2^k-2 \ \text{times}} \\
      \ge & T-(2^k-1). \\
 \end{split}
 \end{equation}
Now, we identify the values of $n$ for which these bounds are attained 
by computing the metric dimension of $\mathcal{E}_{\mathbb{Z}_{n}}$. 
 
\begin{theorem}
   Let  $n= p_{1}^{m_1}p_{2}^{m_2} \cdots p_{k}^{m_k}$, where $p_{1}<p_{2}< \cdots < p_{k}$ are primes, $k\ge 2$, and $m_i>1$ for at least one $i$. Then   \begin{equation*}dim (\mathcal{E}_{\mathbb{Z}_{n}})=
       \begin{cases}
           T-(2^k-1), & \text{if}\ m_i>1\ \text{for at least two}\ i,\\
           T-(2^k-2), & \text{if}\ m_i>1\ \text{for exactly one}\ i.\\
       \end{cases}
   \end{equation*}
\end{theorem}
\begin{proof}
    By Equation $(1)$, we see that any resolving set of $\mathcal{E}_{\mathbb{Z}_{n}}$ must contain at least $T-(2^k-1)$ vertices consisting of all but at most one vertex of each of the distance similar partitioned sets. 
    Case $1$:  $m_i>1$ for at least two $i$\\
    Here, it remains to show that there exists a resolving set of cardinality $T-(2^k-1)$. Take $W$ as an ordered set consisting of $m-1$ vertices of $X$, followed by $|X_i|-1$ vertices of the sets $X_i$, for $1\le i \le 2^k-2$. Without loss of generality, let 
    \begin{equation*} \begin{split}
     W= & X\backslash\{\langle p_{1}^{m_1-1}p_{2}^{m_2-1} \cdots p_{k}^{m_k-1}\rangle\}\bigcup X_1\backslash\{\langle p_{1}^{m_1}p_{2}^{m_2-1} \cdots p_{k}^{m_k-1}\rangle \}\bigcup \cdots \\
      & \bigcup X_{2^k-2}\backslash\{\langle p_{1}^{m_1-1}p_{2}^{m_2} \cdots p_{k}^{m_k}\rangle\}.\\ 
    \end{split}
    \end{equation*}Since $\langle p_{1}^{m_1-1}p_{2}^{m_2-1} \cdots p_{k}^{m_k-1}\rangle$ is the only essential ideal of the set $V\backslash W$,
  we see that  
    \begin{equation*}
       r(\langle p_{1}^{m_1-1}p_{2}^{m_2-1} \cdots p_{k}^{m_k-1}\rangle|W)= (1,1,\cdots, 1)\ne r(v|W)\ \text{for any}\ v\in V\backslash W. 
    \end{equation*} 
    Now, take any two vertices $\hat{I}$ and $\hat{J}$ of $V\backslash W$ with respective index sets $\Xi_{\hat{I}}$ and  $\Xi_{\hat{J}}$.
    we claim that $r(\hat{I}|W)\ne r(\hat{J}|W)$. For $\hat{I}$ and $\hat{J}$, either   $\Xi_{\hat{I}}\cap \Xi_{\hat{J}}= \phi$ or  $\Xi_{\hat{I}}\cap \Xi_{\hat{J}}\ne \phi$. If  $\Xi_{\hat{I}}\cap \Xi_{\hat{J}}= \phi$, then there exist at least two distinct primes $p_i$ and $p_j$ such that $p_i^{m_i}\in \hat{I}\ \text{but} \notin \hat{J}$ and $p_j^{m_j}\in \hat{J}\ \text{but}\ \notin \hat{I}$. Now, $d(\hat{I}, v)=2$
and $d(\hat{J}, v)=1$ for any $v\in X_{\langle p_i^{m_i} \rangle}$. Consequently, $r(\hat{I}|W)$ will have $2$ in all the co-ordinates corresponding to the elements from the set $X_{\langle p_i^{m_i} \rangle}$ whereas $r(\hat{J}|W)$ will have $1$ in the respective coordinates. So $r(\hat{I}|W)\ne r(\hat{J}|W)$. If $\Xi_{\hat{I}}\cap \Xi_{\hat{J}}\ne \phi$, then it can be either $\Xi_{\hat{I}}$ or $\Xi_{\hat{J}}$ or none of these. Let $\Xi_{\hat{I}}\cap \Xi_{\hat{J}}=\Xi_{\hat{I}}$. Since $\Xi_{\hat{J}}\ne \Xi_{\hat{I}}$, there exists at least one prime $p_j$ such that $p_j^{mj}$ is in $\hat{J}$ but not in $\hat{I}$. Hence, $r(\hat{I}|W)$ will have $1$ in all the co-ordinates corresponding to the elements from the set $X_{\langle p_j^{m_j} \rangle}$ and $r(\hat{J}|W)$ will have $2$ in all the co-ordinates corresponding to the elements from the set $X_{\langle p_j^{m_j} \rangle}$. Thus, $r(\hat{I}|W)\ne r(\hat{J}|W)$. Through a similar argument, we see that  $r(\hat{I}|W)\ne r(\hat{J}|W)$ whenever $\Xi_{\hat{I}}\cap \Xi_{\hat{J}}=\Xi_{\hat{J}}$. Now, in the last case, there must exist at least two primes $p_i$ and $p_j$ such that $p_i^{m_i}$ is in $\hat{I}$ but not in $\hat{J}$ and $p_j^{m_j}$ is in $\hat{J}$ but not in $\hat{I}$ leading to $r(\hat{I}|W)\ne r(\hat{J}|W)$. \\
Case $2$: If $m_i>1$ for exactly one $i$\\
 without loss of generality take $n= p_{1}^{m_1}p_2\cdots p_k$, $m_1>1$. We know that any resolving set contains all but at most one vertex of each of the distance similar partitioned sets and by Equation $(1)$, $|W|\ge T-(2^k-1)$, for any resolving set $W$. At first, we show that there is no resolving set of cardinality $T-(2^k-1)$. For this, take $W$ as an ordered set consisting of $m-1$ vertices of $X$ followed by $|X_i|-1$ vertices of $X_i$ for each $i$. That is, 
  \begin{equation*}
  \begin{split}
  W= & \{\langle p_1 \rangle, \langle p_1^2 \rangle, \langle p_1^{m_1-2} \rangle, \langle p_2 \rangle, \langle p_1p_2 \rangle, \cdots, \langle p_1^{m_1-2}p_2 \rangle,  \cdots,  \langle p_k\rangle, \langle p_1p_k \rangle, \cdots, \langle p_1^{m_1-2}p_k \rangle, \\ 
  &
  \langle p_2p_3 \rangle,  \langle p_1p_2p_3 \rangle, \cdots, \langle p_1^{m_1-2} p_2p_3 \rangle, \cdots, \langle p_2p_3\cdots p_k \rangle, \cdots, \langle p_1^{m_1-2}p_2p_3\cdots p_k\rangle \}. \\ 
  \end{split}
  \end{equation*}
Consider the vertices $\hat{I}= \langle p_1^{m_1-1} \rangle$ and $\hat{J}= \langle p_1^{m_1} \rangle$ of $V\backslash W$. Since $\hat{I}$ is essential, $r(\hat{I}|W)= (1,1,\cdots,1)$. For $\hat{J},\  \Xi_{\hat{J}}=\{1\}$ and for any vertex $w\in W$, $1\notin \Xi_{w}$. Consequently, $d(\hat{J}, w)= 1$ and $r(\hat{J}|W)= (1,1,\cdots,1)= r(\hat{I}|W)$. Thus, there is no resolving set of cardinality  $T-(2^k-1)$. Now, take $W'$ as an ordered set obtained by adjoining one more vertex (say $\langle p_1^{m_1}\rangle$) to $W$.
Let  $\hat{I}$ and $\hat{J}$ be two distinct vertices of $V\backslash W$ with index sets $\Xi_{\hat{I}}$ and  $\Xi_{\hat{J}}$ respectively. Then there may occur two cases- either $\Xi_{\hat{I}}\cap \Xi_{\hat{J}}= \phi$ or  $\Xi_{\hat{I}}\cap \Xi_{\hat{J}}\ne \phi$. Proceeding in the same manner as in the proof of case $1$, we see that $W'$ is a resolving set of  $\mathcal{E}_{\mathbb{Z}_{n}}$ of minimum cardinality. 
\end{proof}
\begin{corollary}
    Let $n= p_1^{m_{1}}p_2^{m_{2}}$, where $p_1<p_2$ be primes. Then 
\begin{equation*}
dim(\mathcal{E}_{\mathbb{Z}_{n}})=
\begin{cases}
    2m-2, & \textit{if} \  m_{1}=m \ge 2,  m_{2}= 1 \textit{or vice versa}, \\
    m_{1}m_{2}+m_{1}+m_{2}-4, & \textit{if} \ m_{1}, m_{2}> 1.
\end{cases}.
\end{equation*} 
\end{corollary}

\begin{example}
\begin{itemize}
    \item Consider the graph $\mathcal{E}_{\mathbb{Z}_{n}}$ for $n= p_1^{m_1}p_2^{m_2}p_3$, $m_1,m_2>1$. Then the distance similar partition of vertices is given by,
\begin{equation*}
\begin{split}
    X= & \{\langle p_1^{r_1}p_2^{r_2} \rangle: 0\le i \le m_i-1\ \text{for}\ i=1,2\}\backslash \mathbb{Z}_{n}, |X|=m_1m_2-1\\
    X_1=&\{ \langle p_1^{m_1}p_2^{r_2} \rangle: 0\le r_2 \le m_2-1\},  |X_1|=m_2, \\
    X_2= &\{ \langle p_1^{r_1}p_2^{m_2} \rangle: 0\le r_1 \le m_1-1\}, |X_2|=m_1,  \\
   X_3= & \{ \langle p_1^{r_1}p_2^{r_2}p_3 \rangle: 0\le r_i \le m_i-1\ \text{for}\ i=1,2\}, |X_3|=m_1m_2, \\      
    X_4= &\{ \langle p_1^{m_1}p_2^{m_2} \rangle\},\\
   X_5=& \{ \langle p_1^{m_1}p_2^{r_2}p_3 \rangle: 0\le r_2 \le m_2-1\}, |X_5|=m_2,  \\
    X_6= & \{ \langle p_1^{r_1}p_2^{m_2}p_3 \rangle: 0\le r_1 \le m_1-1\}, |X_6|=m_1.  \\
\end{split}
\end{equation*}
Since any resolving set $W$ contains all but at most one vertex of each of the distance similar vertex partitioned sets,  take $W$ as follows:
\begin{equation*}
\begin{split}
    W= & X\backslash \{\langle p_1^{m_1-1}p_2^{m_2-1} \rangle \}\bigcup X_1\backslash \{\langle p_1^{m_1}p_2^{m_2-1}\rangle\}\bigcup  X_2\backslash \{\langle p_1^{m_1-1}p_2^{m_2}\rangle\}\bigcup   X_3\backslash \{\langle p_1^{m_1-1}p_2^{m_2-1}p_3 \rangle\}\\
    & \bigcup  X_5\backslash \{\langle p_1^{m_1}p_2^{m_2-1}p_3 \rangle\}\bigcup X_6\backslash \{\langle p_1^{m_1-1}p_2^{m_2}p_3 \rangle\}, |W|= 2(m_1m_2+m_1+m_2)-7= T-7\\      
\end{split}
\end{equation*}
To prove $W$ is a  minimum resolving set, it is enough to show that each vertex of  $V\backslash W$ has a unique metric representation.
The representations of the seven vertices of $V\backslash W$ are given as follows:
\begin{equation*}
\begin{split}
r(\langle p_1^{m_1-1}p_2^{m_2-1} \rangle |W) = & (\underbrace{1,1,\cdots,1}_{T-7 \text{times}}),\\
r(\langle p_1^{m_1}p_2^{m_2-1} \rangle |W) = & (\underbrace{1,1,\cdots,1}_{m_1m_2-2}, \underbrace{2,2,\cdots,2}_{m_2-1},\underbrace{1,1,\cdots,1}_{m_1-1},\underbrace{1,1,\cdots,1}_{m_1m_2-1},\underbrace{2,2,\cdots,2}_{m_2-1},\underbrace{1,1,\cdots,1}_{m_1-1}),\\
r(\langle p_1^{m_1-1}p_2^{m_2} \rangle |W) = & (\underbrace{1,1,\cdots,1}_{m_1m_2-2}, \underbrace{1,1,\cdots,1}_{m_2-1},\underbrace{2,2,\cdots,2}_{m_1-1},\underbrace{1,1,\cdots,1}_{m_1m_2-1},\underbrace{1,1,\cdots,1}_{m_2-1},\underbrace{2,2,\cdots,2}_{m_1-1}),\\
r(\langle p_1^{m_1-1}p_2^{m_2-1}p_3 \rangle |W) = & (\underbrace{1,1,\cdots,1}_{m_1m_2-2}, \underbrace{1,1,\cdots,1}_{m_2-1},\underbrace{1,1,\cdots,1}_{m_1-1},\underbrace{2,2,\cdots,2}_{m_1m_2-1},\underbrace{2,2,\cdots,2}_{m_2-1},\underbrace{2,2,\cdots,2}_{m_1-1}),\\
\end{split}
\end{equation*}

\begin{equation*}
\begin{split}
r(\langle p_1^{m_1}p_2^{m_2}\rangle |W) = & (\underbrace{1,1,\cdots,1}_{m_1m_2-2}, \underbrace{2,2,\cdots,2}_{m_2-1},\underbrace{2,2,\cdots,2}_{m_1-1},\underbrace{1,1,\cdots,1}_{m_1m_2-1},\underbrace{2,2,\cdots,2}_{m_2-1},\underbrace{2,2,\cdots,2}_{m_1-1}),\\
r(\langle p_1^{m_1}p_2^{m_2-1}p_3 \rangle |W) = & (\underbrace{1,1,\cdots,1}_{m_1m_2-2}, \underbrace{2,2,\cdots,2}_{m_2-1}, \underbrace{1,1,\cdots,1}_{m_1-1},\underbrace{2,2,\cdots,1}_{m_1m_2-1},\underbrace{2,2,\cdots,2}_{m_2-1},\underbrace{2,2,\cdots,2}_{m_1-1}),\\
r(\langle p_1^{m_1-1}p_2^{m_2}p_3 \rangle |W) = & (\underbrace{1,1,\cdots,1}_{m_1m_2-2}, \underbrace{1,1,\cdots,1}_{m_2-1}, \underbrace{2,2,\cdots,2}_{m_1-1},\underbrace{2,2,\cdots,1}_{m_1m_2-1},\underbrace{2,2,\cdots,2}_{m_2-1},\underbrace{2,2,\cdots,2}_{m_1-1}).\\
\end{split}
\end{equation*}
It can be seen that any two distinct vertices of $V(\mathcal{E}_{\mathbb{Z}_{n}})\backslash W$ have different metric representations with respect to $W$. Thus $W$ is a resolving set having $T-7=2(m_1m_2+m_1+m_2)-7$ vertices. Also, any resolving set must contain more than $T-7$ elements, we conclude that $dim (\mathcal{E}_{\mathbb{Z}_{n}})= T-7$.
\item Let $n= p_1^{m_1}p_2p_3$, $m_1>1$. Then the distance similar partition of vertices of $\mathcal{E}_{\mathbb{Z}_{n}}$ is given by:
\begin{equation*}
\begin{split}
    X= & \{\langle p_1^{r_1}\rangle: 1\le r_1 \le m_1-1\}, |X|=m_1-1,\\
    X_1=&\{ \langle p_1^{m_1} \rangle\},\\
    X_2= &\{ \langle p_1^{r_1}p_2 \rangle: 0\le r_1 \le m_1-1\}, |X_2|=m_1,  \\
   X_3=& \{ \langle p_1^{r_1}p_3 \rangle: 0\le r_1 \le m_1-1\}, |X_3|=m_1,  \\      
    X_4= &\{ \langle p_1^{m_1}p_2\rangle\}, \ 
   X_5= \{ \langle p_1^{m_1}p_3 \rangle\},\\
    X_6= & \{ \langle p_1^{r_1}p_2p_3 \rangle: 0\le r_1 \le m_1-1\}, |X_6|=m_1.  \\
\end{split}
\end{equation*}
Now, take $W$as an ordered set consisting of first $m_1-2$ vertices of $X$, first $m_1-1$ vertices of $X_2$, and so on.
Then, for the vertices $\langle p_1^{m_1-1}\rangle$ and $\langle p_1^{m_1}\rangle$  of $V(\mathcal{E}_{\mathbb{Z}_{n}})\backslash W$, 
\begin{equation*}
r(\langle p_1^{m_1-1}\rangle|W)= (\underbrace{1,1,\cdots,1}_{m_1-2}, \underbrace{1,1,\cdots,1}_{m_1-1},\underbrace{1,1,\cdots,1}_{m_1-1},\underbrace{1,1,\cdots,1}_{m_1-1})= r(\langle p_1^{m_1}\rangle|W)    
\end{equation*}
Hence, $W$ cannot be a resolving set of $\mathcal{E}_{\mathbb{Z}_{n}}$. Without loss of generality, take $W'=W\cup \{\langle p_1^{m_1}\rangle \}$. That is, \[ W'= X\backslash\{ \langle p_1^{m_1-1}\rangle\}\bigcup X_1\bigcup X_2\backslash\{ \langle p_1^{m_1-1}p_2 \rangle\}\bigcup X_3\backslash\{ \langle p_1^{m_1-1}p_3 \rangle\}\bigcup X_6\backslash\{ \langle p_1^{m_1-1}p_2p_3 \rangle\}.\] The representations of vertices of $V(\mathcal{E}_{\mathbb{Z}_{n}})\backslash W'$ are given by,
\begin{equation*}
\begin{split}
  r(\langle p_1^{m_1-1}\rangle|W')= & (\underbrace{1,1,\cdots,1}_{m_1-2}, 1,\underbrace{1,1,\cdots,1}_{m_1-1},\underbrace{1,1,\cdots,1}_{m_1-1},\underbrace{1,1,\cdots,1}_{m_1-1}), \\
  r(\langle p_1^{m_1-1}p_2\rangle|W')= & (\underbrace{1,1,\cdots,1}_{m_1-2}, 1,\underbrace{2,2,\cdots,2}_{m_1-1},\underbrace{1,1,\cdots,1}_{m_1-1},\underbrace{2,2,\cdots,2}_{m_1-1}), \\
  r(\langle p_1^{m_1-1}p_3\rangle|W')= & (\underbrace{1,1,\cdots,1}_{m_1-2}, 1,\underbrace{1,1,\cdots,1}_{m_1-1},\underbrace{2,2,\cdots,2}_{m_1-1},\underbrace{2,2,\cdots,2}_{m_1-1}), \\ 
  r(\langle p_1^{m_1}p_2\rangle|W')= & (\underbrace{1,1,\cdots,1}_{m_1-2}, 2,\underbrace{2,2,\cdots,2}_{m_1-1},\underbrace{1,1,\cdots,1}_{m_1-1},\underbrace{2,2,\cdots,2}_{m_1-1}), \\ 
 r(\langle p_1^{m_1}p_3\rangle|W')= & (\underbrace{1,1,\cdots,1}_{m_1-2}, 2,\underbrace{1,1,\cdots,1}_{m_1-1},\underbrace{2,2,\cdots,2}_{m_1-1},\underbrace{2,2,\cdots,2}_{m_1-1}), \\ 
r(\langle  p_1^{m_1-1}p_2p_3\rangle|W')= & (\underbrace{1,1,\cdots,1}_{m_1-2}, 1,\underbrace{2,2,\cdots,2}_{m_1-1},\underbrace{2,2,\cdots,2}_{m_1-1},\underbrace{2,2,\cdots,2}_{m_1-1}). \\ 
\end{split}
\end{equation*}
\end{itemize}
This unique representation of vertices of $V(\mathcal{E}_{\mathbb{Z}_{n}})\backslash W'$  ensures that $W'$ is a minimum resolving set of $\mathcal{E}_{\mathbb{Z}_{n}}$. Hence, the  metric dimension of $\mathcal{E}_{\mathbb{Z}_{n}}$ is $4(m_1-1)=T-6$.
\end{example}

\section{Zagreb Indices of $\mathcal{E}_{\mathbb{Z}_{n}}$}
In this section, we calculate the $1$st and $2$nd Zagreb indices of $\mathcal{E}_{\mathbb{Z}_{n}}$.\\
\begin{proposition}
    Let $n=p^m$, $m>2$. Then \begin{enumerate}
        \item The first Zagreb index of $\mathcal{E}_{\mathbb{Z}_{n}}$ is, $M_1(\mathcal{E}_{\mathbb{Z}_{n}})= (m-1)(T-1)^2$.
        \item The second Zagreb index of $\mathcal{E}_{\mathbb{Z}_{n}}$ is, $M_2(\mathcal{E}_{\mathbb{Z}_{n}})= \binom{m-1}{2}(T-1)^2$.
    \end{enumerate}
\end{proposition}
\begin{lemma}\label{adjin prdct of primes}\cite{jamsheena2023adjacency}
Let $n= p_{1}p_{2} \cdots p_{k}$, where $p_{1}, p_{2}, \cdots, p_{k}$ are distinct primes. Then any two vertices $\langle x \rangle$ and $\langle y \rangle$ of the essential ideal graph of $\mathbb{Z}_{n}$ are adjacent if and only if $gcd(x,y)=1$, provided $x$ is the product of $i$ distinct primes and $y$ is the product of $j$ distinct primes for $1 \le i,j \le k-1 $.
\end{lemma}

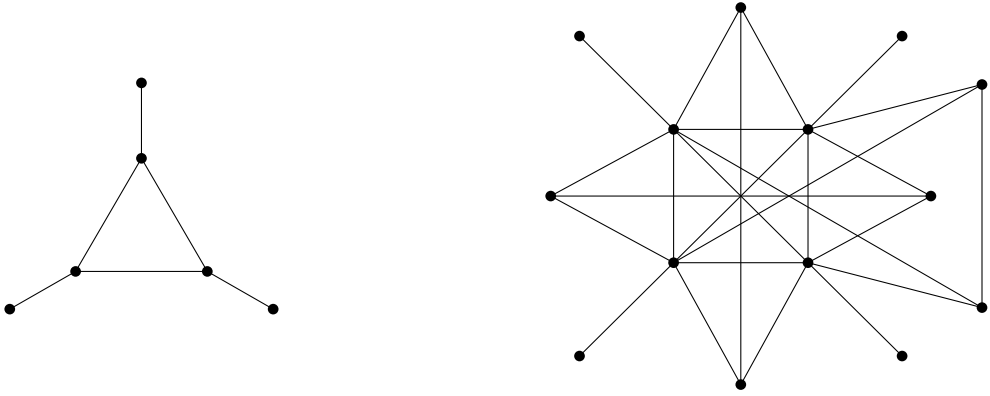
\begin{figure}
		\begin{minipage}{7 cm}
			\begin{tikzpicture}[scale=0.5]\label{fig3.1}
			\fill [black] (0,0) +(90:2) circle (4pt);
			\fill [black](0,0) +(210:2) circle (4pt);
			\fill [black] (0,0) +(330:2) circle (4pt);
			\fill [black] (0,0) +(90:4) circle (4pt);
			\fill [black] (0,0) +(210:4) circle (4pt);
			\fill [black] (0,0) +(330:4) circle (4pt);
			\draw (90:4)--(90:2)--(210:2)--(330:2)--(90:2);
			\draw (210:2)--(210:4);
			\draw (330:2)--(330:4);
			\end{tikzpicture}
		\end{minipage}
		\begin{minipage}{5.8 cm}
			\begin{tikzpicture}[scale=0.5]
			\fill [black] (0,0) +(45:2.5) circle (4pt);
			\fill [black] (0,0) +(135:2.5) circle (4pt);
			\fill [black] (0,0) +(225:2.5) circle (4pt);
			\fill [black] (0,0) +(315:2.5) circle (4pt);
			\fill [black] (0,0) +(45:6) circle (4pt);
			\fill [black] (0,0) +(135:6) circle (4pt);
			\fill [black] (0,0) +(225:6) circle (4pt);
			\fill [black] (0,0) +(315:6) circle (4pt);
			\fill [black] (0,0) +(90:5) circle (4pt);
			\fill [black] (0,0) +(180:5) circle (4pt);
			\fill [black] (0,0) +(270:5) circle (4pt);
			\fill [black] (0,0) +(0:5) circle (4pt);
			\draw (45:2.5)--(135:2.5)--(225:2.5)--(315:2.5)--(45:2.5)--(45:6);
			\draw (135:2.5)--(135:6);
			\draw (225:2.5)--(225:6);
			\draw (315:2.5)--(315:6);
			\draw (135:2.5)--(315:2.5);
			\draw (225:2.5)--(45:2.5);
			\draw (135:2.5)--(90:5)--(45:2.5)--(0:5)--(315:2.5)--(270:5)--(225:2.5)--(180:5)--(135:2.5);
			\draw (270:5)--(90:5);
			\draw (0:5)--(180:5);
			\fill [black] (0,0) +(25:7) circle (4pt);
			\fill [black] (0,0) +(-25:7) circle (4pt);
			\draw (25:7)--(-25:7);
			\draw (45:2.5)--(25:7)--(225:2.5);
			\draw(135:2.5)--(-25:7)--(315:2.5);
			\end{tikzpicture}
		\end{minipage}
		\caption{$\mathcal{E}_{\mathbb{Z}_{n}}$ for $n= \prod_{i=1}^3p_i$ and $n= \prod_{i=1}^4p_i$} \label{Figure 1}
			\end{figure}

\begin{theorem}\label{zagreb of distct primes}
    Let $n=p_1p_2\cdots p_k$. Then,  the fist and second Zagreb indices of $\mathcal{E}_{\mathbb{Z}_{n}}$ are, \begin{enumerate}
        \item $M_1(\mathcal{E}_{\mathbb{Z}_{n}})= \displaystyle\sum_{i=1}^{k-1}\binom{k}{i} (2^{k-i}-1)^2.$
        \item $M_2(\mathcal{E}_{\mathbb{Z}_{n}})= \displaystyle\sum_{t=1}^{\lfloor \frac{k}{2} \rfloor}\binom{k}{t}(2^{k-t}-1) [\frac{1}{2} \binom{k-t}{t}(2^{k-t}-1)+ \displaystyle\sum_{s=t}^{k-t}\binom{k-t}{s}(2^{k-s}-1)].$
    \end{enumerate}
\end{theorem}
\begin{proof}
     For $n= p_1p_2\cdots p_k$, the vertex set of $\mathcal{E}_{\mathbb{Z}_{n}}$ can be partitioned as follows: \\
\begin{equation*}
    \begin{split}
      & V_1=\{\langle p_i\rangle : 1\le i \le k \}
\\
& V_2=\{\langle p_ip_j\rangle: 1 \le i \le k-1 \textit{and} \  i+1 \le j \le k \}
\\
& V_3=\{\langle p_ip_jp_l\rangle : 1 \le i \le k-2, \ i+1 \le j \le k-1 \textit{and} \  j+1 \le l \le k \}
\\
  & \vdots
  \\
& V_{k-1}=\{\langle p_1p_2p_3\cdots p_{k-1}\rangle, \langle p_1p_2p_3\cdots p_{k-2}p_k\rangle,\cdots, \langle p_2p_3\cdots p_{k-1}p_k\rangle \}
    \end{split}
\end{equation*}

Clearly $\vert V_1 \vert = \binom{k}{1}, \ |V_2|= \binom{k}{2}, \cdots, \ \text{and} \  |V_{k-1}|= \binom{k}{k-1}$. Also, by Lemma \ref{adjin prdct of primes}, two vertices $\langle x \rangle$ and $\langle y \rangle$ of $\mathcal{E}_{\mathbb{Z}_{n}}$ are adjacent if and only if the generators $x$ and $y$ have no prime factors in common.
\begin{equation*} \text{For a vertex}\ v \in V(\mathcal{E}_{\mathbb{Z}_{n}}),\
deg(v) = \begin{cases}  
2^{k-1}-1, & \textit{if}\ v\in V_1, \\
2^{k-2}-1, & \textit{if} \ v\in V_2, \\
\vdots    &    \vdots \\
3, &  \textit{if} \  v\in V_{k-2}, \\
1, &  \textit{if} \  v\in V_{k-1}. \\
\end{cases}
\end{equation*}
Also, for a fixed $i$, 
 \[\sum_{v\in V_i} deg(v)^2= \binom{k}{i} (2^{k-i}-1)^2.\] Hence,
 \begin{enumerate}
    \item \begin{equation*}\begin{split}
     M_1(\mathcal{E}_{\mathbb{Z}_{n}})= & \sum_{v\in V_1}deg(v)^2+ \sum_{v\in V_2}deg(v)^2+\cdots+  \sum_{v\in V_{k-1}}deg(v)^2\\
       = &  \sum_{i=1}^{k-1} \sum_{v\in V_i}deg(v)^2 \\
     = &  \sum_{i=1}^{k-1}\binom{k}{i} (2^{k-i}-1)^2.\\
\end{split}
\end{equation*}
 \item  For a fixed $t$, $1\le t\le \lfloor \frac{k}{2} \rfloor$,  Lemma \ref{adjin prdct of primes} assures that   any vertex  $u \in V_t$ are adjacent only to $\binom{k-t}{s}$ vertices of $V_s$ for $t\le s\le k-t$.
Then, for $u,v \in V_t, \ 1\le t\le \lfloor \frac{k}{2} \rfloor$, 
\begin{equation}\label{u,v in Vt}
    \displaystyle\sum_{u\sim v} deg(u)deg(v)= \frac{1}{2}\binom{k}{t}\binom{k-t}{t}(2^{k-t}-1)^2.
    \end{equation}  
Now, consider $u\in V_t$ for a fixed $t$ such that $1\le t\le \lfloor \frac{k}{2} \rfloor$. Then we have,
\begin{equation}\label{v in Vs}
\begin{split}
    \displaystyle\sum_{u\sim v \atop v\in V_s}deg(u)deg(v)= &\displaystyle\sum_{u\sim v \atop v\in V_t}deg(u)deg(v)+  \displaystyle\sum_{u\sim v \atop v\in V_{t+1}}deg(u)deg(v)+ \cdots+
   \displaystyle\sum_{u\sim v \atop v\in V_{k-t}}deg(u)deg(v)\\ 
    = & \frac{1}{2}\binom{k}{t}\binom{k-t}{t}(2^{k-t}-1)^2+ \binom{k}{t}\binom{k-t}{t+1}(2^{k-t}-1)(2^{k-(t+1)}-1)+ \cdots\\
    + & \binom{k}{t}\binom{k-t}{k-t}(2^{k-t}-1)(2^{t}-1)\\
    = & \frac{1}{2}\binom{k}{t}\binom{k-t}{t}(2^{k-t}-1)^2+ \binom{k}{t} (2^{k-t}-1)\displaystyle\sum_{s=t+1}^{k-t}(2^{k-s}-1).
\end{split}
\end{equation}
Now, by Equation
(\ref{v in Vs}),
\begin{equation*} 
\begin{split}
    M_2(\mathcal{E}_{\mathbb{Z}_{n}})= &   \displaystyle\sum_{u\sim v\atop u,v\in V(\mathcal{E}_{\mathbb{Z}_{n}})}deg(u)deg(v) \\
    = & \displaystyle\sum_{u\sim v\atop u\in V_1, \ v\in V_s}deg(u)deg(v)+ \displaystyle\sum_{u\sim v\atop u\in V_2, \ v\in V_s}deg(u)deg(v)+\cdots+ \displaystyle\sum_{u\sim v\atop u\in V_{\lfloor \frac{k}{2} \rfloor}, \ v\in V_s} deg(u)deg(v)\\
   = & \frac{1}{2} \displaystyle\sum_{t=1}^{\lfloor \frac{k}{2} \rfloor}\binom{k}{t}\binom{k-t}{t}(2^{k-t}-1)^2+  \displaystyle\sum_{t=1}^{\lfloor \frac{k}{2} \rfloor}\binom{k}{t}(2^{k-t}-1) \displaystyle\sum_{s=t}^{k-t}\binom{k-t}{s}(2^{k-s}-1)\\
   =& \displaystyle\sum_{t=1}^{\lfloor \frac{k}{2} \rfloor}\binom{k}{t}(2^{k-t}-1) [\frac{1}{2} \binom{k-t}{t}(2^{k-t}-1)+ \displaystyle\sum_{s=t}^{k-t}\binom{k-t}{s}(2^{k-s}-1)].
\end{split}
\end{equation*}
 \end{enumerate}
\end{proof} 
\begin{example}
\begin{itemize}
      \item Let $n=p_1p_2p_3$. Then by Theorem \ref{zagreb of distct primes}, the first and second Zagreb indices are given by,
      \begin{equation*}\begin{split}
     M_1(\mathcal{E}_{\mathbb{Z}_{n}})= &  \displaystyle\sum_{i=1}^{2}\binom{3}{i} (2^{3-i}-1)^2  \\       
    = & \binom{3}{1} 3^2+ \binom{3}{2} = 30, 
    \end{split}
    \end{equation*} and
     \begin{equation*}\begin{split}
     M_2(\mathcal{E}_{\mathbb{Z}_{n}})= & \displaystyle\sum_{t=1}^{\lfloor \frac{3}{2} \rfloor}\binom{3}{t}(2^{3-t}-1)\displaystyle\sum_{s=t}^{3-t}\binom{3-t}{s}(2^{3-s}-1) \\
     = & \binom{3}{1} 3[3 \binom{2}{1}+1]= 63.
     \end{split}
    \end{equation*} 
   
    \item Let $n=p_1p_2p_3p_4$. Then 
    \begin{equation*}
    \begin{split}
     M_1(\mathcal{E}_{\mathbb{Z}_{n}})= &  \displaystyle\sum_{i=1}^{3}\binom{4}{i} (2^{4-i}-1)^2  \\       
    = & \binom{4}{1} 7^2 +\binom{4}{2} 3^2+ \binom{4}{3} = 254, 
    \end{split}
    \end{equation*} and 
    \begin{equation*}\begin{split}
     M_2(\mathcal{E}_{\mathbb{Z}_{n}})= & \displaystyle\sum_{t=1}^2 \binom{4}{t}(2^{4-t}-1)\displaystyle\sum_{s=t}^{4-t}\binom{4-t}{s}(2^{4-s}-1) \\
     = & \binom{4}{1} 7[ \binom{3}{1} 7+3\binom{3}{2}+ 1]+  \binom{4}{2} 9= 922.
     \end{split}
    \end{equation*} 
\end{itemize} 
\end{example}
 This can be easily verified from Figure \ref{Figure 1}.

\begin{theorem}\label{Zagreb for any n}
    Let  $n= p_{1}^{m_1}p_{2}^{m_2} \cdots p_{k}^{m_k}$, where $p_{1}<p_{2}< \cdots < p_{k}$ are primes, $k\ge 2$ and  $m_i>1$ for at least one $i$. Then, \begin{enumerate}
        \item $M_1(\mathcal{E}_{\mathbb{Z}_{n}})=  m(T-1)^2+ \displaystyle\sum_{\hat{I}}|X_{\hat{I}}|(m+\displaystyle \sum_{\hat{J}: \ \Xi_{\hat{I}}\cap \Xi_{\hat{J}}=\phi}|X_{\hat{J}}|)^2$.
         \item  \begin{equation*}   \begin{split}
              M_2(\mathcal{E}_{\mathbb{Z}_{n}})=  & \binom{m}{2}(T-1)^2 + m(T-1)  \displaystyle\sum_{\hat{I}} |X_{\hat{I}}|(m+\displaystyle\sum_{\hat{J}:\ \Xi_{\hat{I}}\cap \Xi_{\hat{J}}=\phi}|X_{\hat{J}}|)\\
      + & \frac{\delta_i^j}{2} |X_{\hat{I}}||X_{\hat{J}}|\displaystyle\sum_{u,v\in \mathscr{U}} deg(u)deg(v), \text{where}\quad \delta_i^j=\begin{cases}
        1, & \text{if}\ \Xi_{\hat{I}}\cap \Xi_{\hat{J}}=\phi,\\
        0, & \text{otherwise},
    \end{cases} \end{split}
\end{equation*}and $deg(u)= |X|+\displaystyle \sum_{\hat{J}: \ \Xi_{\hat{I}}\cap \Xi_{\hat{J}}=\phi}|X_{\hat{J}}|$ for any vertex $u\in X_{\hat{I}}$ of $\mathscr{U}$.
\end{enumerate}
\end{theorem}
\begin{proof}
   By Theorem \ref{Complt Strctr of EZn}, $\mathcal{E}_{\mathbb{Z}_{n}}\cong K_m \vee \mathscr{G}[\Gamma_1,\Gamma_2,\cdots,\Gamma_{2^k-2}]$, 
 where $K_m$ is the subgraph induced by the set $X$ of essential ideals of $\mathbb{Z}_{n}$, and  $\Gamma_i= \mathcal{E}_{\mathbb{Z}_{n}}(X_{\hat{I}})$ for each of the  $2^k-2$ equivalent class $ X_{\hat{I}}$ of the set $\mathscr{U}$.
 Thus, 
 $V(\mathcal{E}_{\mathbb{Z}_{n}})=  X \cup _{\hat{I}}X_{\hat{I}}$, where the union is taken over all the equivalent classes.
 Then, for any vertex $v\in X$, $deg(v)=T-1$ and by Lemma \ref{Chrczn Adcncy}, 
 \begin{equation}\label{deg in Xi}
\text{for any vertex}\ v\in X_{\hat{I}}, \hspace{1.5cm} deg(v)= |X|+\displaystyle \sum_{\hat{J}: \ \Xi_{\hat{I}}\cap \Xi_{\hat{J}}=\phi}|X_{\hat{J}}|,\end{equation}
 \begin{enumerate}
\item \begin{equation*}
     \begin{split}
         M_1(\mathcal{E}_{\mathbb{Z}_{n}})= & \displaystyle \sum_{v \in X}deg(v)^2+ \displaystyle \sum_{v \in \cup_{\hat{I}}X_{\hat{I}}}  deg(v)^2, \\
       = & m(T-1)^2+ \displaystyle\sum_{\hat{I}}|X_{\hat{I}}|(m+\displaystyle \sum_{\hat{J}: \ \Xi_{\hat{I}}\cap \Xi_{\hat{J}}=\phi}|X_{\hat{J}}|)^2,
     \end{split}
 \end{equation*}
 where the summation runs over all the equivalent classes.
\item  First consider $v\in X$. Since, $|X|=m$ and there are $\binom{m}{2}$ pairs of elements in $X$ , we have 
\begin{equation}\label{deg 1}
    \displaystyle\sum_{u,v\in X} deg(u)deg(v)= \binom{m}{2}(T-1)^2.  
\end{equation}
Now, each vertex $u \in X$ is adjacent to every vertex  $v\in X_{\hat{I}}$, for each of the equivalent class $X_{\hat{I}}$. Then, 
\begin{equation*}
    deg(u)deg(v)= (T-1)(m+\displaystyle \sum_{\hat{J}: \ \Xi_{\hat{I}}\cap \Xi_{\hat{J}}=\phi}|X_{\hat{J}}|),
\end{equation*}   and hence, 
\begin{equation}\label{deg 2}
\begin{split}
  \displaystyle\sum_{u\in X\atop v\in \mathscr{U}}  deg(u)deg(v)= & \displaystyle\sum_{\hat{I}}m(T-1)|X_{\hat{I}}|(m+\displaystyle\sum_{\hat{J}:\ \Xi_{\hat{I}}\cap \Xi_{\hat{J}}=\phi}|X_{\hat{J}}|)\\
    = &  m(T-1)  \displaystyle\sum_{\hat{I}} |X_{\hat{I}}|(m+\displaystyle\sum_{\hat{J}:\ \Xi_{\hat{I}}\cap \Xi_{\hat{J}}=\phi}|X_{\hat{J}}|).\\
\end{split}
\end{equation}
Now, if we take $u,v$ from the vertex subset $\mathscr{U}=\bigcup\limits_{\hat{I}}X_{\hat{I}}$, then $u\in X_{\hat{I}}$ and $v\in X_{\hat{J}}$ for some equivalent classes $X_{\hat{I}}$ and $X_{\hat{J}}$. By Lemma \ref{Chrczn Adcncy}, $u$ and $v$ are adjacent if and only if  $\Xi_{\hat{I}}\cap \Xi_{\hat{J}}= \phi$ and hence 
\begin{equation}\label{deg 3}
     \displaystyle\sum_{u\sim v\atop u,v\in \mathscr{U}}deg(u)deg(v)= \delta_i^j
    |X_{\hat{I}}||X_{\hat{J}}|\displaystyle\sum_{u,v\in \mathscr{U}}deg(u)deg(v), \text{where}\quad \delta_i^j=\begin{cases}
        1, & \text{if}\ \Xi_{\hat{I}}\cap \Xi_{\hat{J}}=\phi,\\
        0, & \text{otherwise}.
    \end{cases}
\end{equation} Using Equations (\ref{deg 1}), (\ref{deg 2}) and (\ref{deg 3}), we have
\begin{equation*} 
    \begin{split}
      M_2(\mathcal{E}_{\mathbb{Z}_{n}})= &   \displaystyle\sum_{u,v\in X} deg(u)deg(v)+  \displaystyle\sum_{u\in X\atop v\in \mathscr{U}}  deg(u)deg(v)+   \displaystyle\sum_{u\sim v\atop u,v\in \mathscr{U}}deg(u)deg(v)\\
      = & \binom{m}{2}(T-1)^2 + m(T-1)  \displaystyle\sum_{\hat{I}} |X_{\hat{I}}|(m+\displaystyle\sum_{\hat{J}:\Xi_{\hat{I}}\cap \Xi_{\hat{J}}=\phi}|X_{\hat{J}}|)\\
      + & \delta_i^j |X_{\hat{I}}||X_{\hat{J}}|\displaystyle\sum_{u,v \in \mathscr{U}}deg(u)deg(v), \text{where}\quad \delta_i^j=\begin{cases}
        1, & \text{if}\ \Xi_{\hat{I}}\cap \Xi_{\hat{J}}=\phi,\\
        0, & \text{otherwise},
    \end{cases} \end{split}
\end{equation*} $deg(u)$ and $deg(v)$ are given by Equation (\ref{deg in Xi}). 
 \end{enumerate}
\end{proof}
\begin{corollary}
Let $n= p_1^{m_{1}}p_2^{m_{2}}$, where  $p_1<p_2$ are primes and $m_i>1$ for at least one $i$. Then
\begin{enumerate}
    \item $M_1(\mathcal{E}_{\mathbb{Z}_{n}})= m(T-1)^2+m_1(m+m_2)^2+m_2(m+m_1)^2$
     \item $M_2(\mathcal{E}_{\mathbb{Z}_{n}})= \binom{m}{2}(T-1)^2+ m(T-1)[m(m_1+m_2)+2m_1m_2]+m_1m_2(m+m_1)(m+m_2).$
\end{enumerate}
\end{corollary}
\begin{proof}
    By Theorem \ref{Complt Strctr of EZn}, $\mathcal{E}_{\mathbb{Z}_{n}}\cong K_m\vee K_2[\overline{K_{m_2}}, \ \overline{K_{m_1}}]$, where $K_m$ is the subgraph induced by the set $X$ of essential ideals of $\mathbb{Z}_{n}$. Here, 
    \begin{equation*} \begin{split}
        X_1=& X_{\langle p_1^{m_1} \rangle}= \{ \langle p_1^{m_1} p_2^{r_2} \rangle : 0\le r_2< m_2\}; |X_1|= m_2,\\
         X_2=&  X_{\langle p_2^{m_2} \rangle}= \{ \langle p_1^{r_1} p_2^{m_2} \rangle : 0\le r_1< m_1\}; |X_2|= m_1. 
    \end{split}
    \end{equation*} Then, by Theorem \ref{Zagreb for any n}, $M_1(\mathcal{E}_{\mathbb{Z}_{n}})= m(T-1)^2+ m_1(m+m_2)^2+m_2(m+m_1)^2$, and
    \begin{equation*}\begin{split}
      M_2(\mathcal{E}_{\mathbb{Z}_{n}})= &    \binom{m}{2}(T-1)^2+ m(T-1)[m_2(m+m_1)+m_1(m+m_2)] +m_1m_2(m+m_1)(m+m_2)\\
       = & \binom{m}{2}(T-1)^2+m(T-1)[m(m_1+m_2)+2m_1m_2]+m_1m_2(m+m_1)(m+m_2).
    \end{split} \end{equation*}
\end{proof}

\begin{example}
    Let $n=p_1^2p_2^3p_3^2$. Then, $T= |V(\mathcal{E}_{\mathbb{Z}_{n}})|= 34$, and $\mathcal{E}_{\mathbb{Z}_{n}}\cong K_m \vee \mathscr{G}[\overline{K_6}, \overline{K_4}, \overline{K_6}, \overline{K_2}, \overline{K_3}, \overline{K_2}]$, where $m=11$. The partitioned sets of nonessential ideals of $\mathbb{Z}_{n}$ are 
    \begin{equation*}
    \begin{split}
        X_1= & X_{\langle p_1^2 \rangle }= \{\langle p_1^2p_2^{r_2}p_3^{r_3}\rangle: 0 \le r_i\le <m_i \ \text{for}\ i=2,3\}; \ |X_1|= 6,\\
        X_2= & X_{\langle p_2^3\rangle }= \{\langle p_1^{r_1}p_2^{3}p_3^{r_3}\rangle: 0 \le r_i< m_i\ \text{for}\ i=1,3\}; \ |X_2|= 4,\\
        X_3= & X_{\langle p_3^2 \rangle }= \{\langle p_1^{r_1}p_2^{r_2}p_3^3 \rangle: 0 \le r_i< m_i\ \text{for}\ i=1,2 \}; \ |X_3|= 6,\\ 
         X_4= & X_{\langle p_1^2p_2^3\rangle }= \{\langle p_1^2p_2^{3}p_3^{r_3}\rangle: 0 \le r_3< m_3 \}; \ |X_4|= 2,\\
          X_5= & X_{\langle p_1^2p_3^2 \rangle }= \{\langle p_1^2p_2^{r_2}p_3^{2}\rangle: 0 \le r_2< m_2 \}; \ |X_5|= 3,\\
         X_6 = & X_{\langle p_2^3p_3^2 \rangle }= \{\langle p_1^{r_1}p_2^3p_3^{2}\rangle: 0 \le r_1< m_1 \}; \ |X_6|= 2.
    \end{split}
    \end{equation*}
    
      $deg(u)= \begin{cases}
       11+ |X_2|+ |X_3|+|X_6|= 23, & \ \text{for}\ u\in X_1, \\
        11+  |X_1|+ |X_3|+ |X_5|= 26, & \ \text{for}\ u\in X_2, \\
        11+  |X_1|+ |X_2|+ |X_4|= 23, & \ \text{for}\ u\in X_3, \\
        11+  |X_3|= 17, & \ \text{for}\ u\in X_4, \\
         11+  |X_2|= 15, & \ \text{for}\ u\in X_5, \\
          11+  |X_1|= 17, & \ \text{for}\ u\in X_6. \\
    \end{cases}$\\
Then, by Theorem, 
\begin{equation*}
\begin{split}
    M_1(\mathcal{E}_{\mathbb{Z}_{n}})= & 11\times 33^2 +6\times23^2 +4 \times26^2 + 6\times 23^2 + 2 \times 17^2 + 3\times 15^2+ 2\times 17^2   \\
     =& 22,862.
\end{split}
\end{equation*}
\begin{equation*}
\begin{split}
    M_2(\mathcal{E}_{\mathbb{Z}_{n}})= & \binom{11}{2}33^2+ 11 \times 33 [6\times 23+ 4\times 26+ 6\times 23+ 2\times17 + 3\times 15+2\times 17]+ \frac{1}{2} [6\times 4\times23\times 26 \\
                                     + & 6\times 6 \times  23\times 23+ 6\times 2\times 23\times 17+ 4\times 6\times 26\times 23+4\times 6\times 26 \times 23+ 4 \times 3\times 26\times 15\\
                                     + &  6\times 6\times 23\times 23+ 6\times 4\times 23\times 26+ 6\times 2\times 23\times 17+  2\times  6\times 17\times 23+ 3\times 4\times15\times 26 \\
                                     + & 2\times 6\times 17\times 23]\\.
                                     =& 3,00,666.
\end{split}
\end{equation*}
\end{example}

\section{Conclusion}
In this article, we have proved that the metric dimension of the essential ideal graph $\mathcal{E}_{R}$ of a commutative ring $R$ is finite whenever each vertex of  $\mathcal{E}_{R}$ is of finite degree. 
Also, for the ring  $\mathbb{Z}_{n}$, it is identified that the  graphs $\mathcal{E}_{\mathbb{Z}_{n}}$ and $\mathbb{AIG}(\mathbb{Z}_{n})$ coincide (up to isomorphism) when $n$ is a product of distinct primes. 
Furthermore, we have calculated the metric dimension of $\mathcal{E}_{\mathbb{Z}_{n}}$. Additionally, an alternative method has been provided to establish an upper limit for $dim (\mathcal{E}_{\mathbb{Z}_{n}})$ when $n= p_1p_2\cdots p_k;\ k\ge 6$. Finally, the first and second Zagreb indices of $\mathcal{E}_{\mathbb{Z}{n}}$ are computed for arbitrary values of $n$.

\section{Declarations}
\textbf{Conflict of interest} On behalf of all authors, the corresponding author states that there is no conflict of interest.

\nocite{*}
\printbibliography
\end{document}